\documentclass[a4paper,12pt]{article}
\usepackage{amsthm}
\usepackage{amssymb}
\usepackage{amsmath}
\usepackage{amsfonts}
\usepackage{color}
\usepackage{graphicx}
\usepackage[square, comma, sort&compress, numbers]{natbib}

\newtheorem{lemma}{Lemma}[section]

\newtheorem{proposition}{Proposition}[section]

\theoremstyle{definition}
\newenvironment{demo*}{\vspace{3mm}\noindent{\bf Proof.}}{\hfill $\Box$ \vspace{3mm}}

\usepackage{amssymb}

\begin{document}

\baselineskip=21pt

\title{\bf \Large {Representations of the Riemann zeta function: A probabilistic approach}}
{\color{red}{\author{\normalsize{Jiamei Liu\;\;\; Yuxia Huang\;\;\; Chuancun Yin}\\
{\normalsize\it  School of Statistics,  Qufu Normal University}\\
\noindent{\normalsize\it Shandong 273165, China}\\
e-mail:  ccyin@qfnu.edu.cn}}}
\maketitle
\vskip0.01cm
\noindent{{\bf Abstract}  In this paper, we give a short elementary proof of the well known  Euler's recurrence formula for the   Riemann zeta
function at  positive even integers and  integral representations of the Riemann zeta function at positive integers and at  fractional points by means of   probabilistic approach.   The proof is based on the  moment generating function and  the characteristic function of  logistic and half-logistic  distributions in probability theory.}




\numberwithin{equation}{section}
\section{Introduction}\label{intro}

The well known Riemann zeta function $\zeta$ is defined by
\begin{eqnarray*}
 \zeta(s)=\left\{\begin{array}{ll} \sum_{n=1}^{\infty}\frac{1}{n^s}=\frac{1}{1-2^{-s}}\sum_{n=1}^{\infty}\frac{1}{(2n-1)^s},  \ & {\rm if}\; {\cal{R}}(s)>1,\\
\frac{1}{1-2^{1-s}}\sum_{n=1}^{\infty}\frac{(-1)^{n+1}}{n^s} ,\ &{\rm if}\;   {\cal{R}}(s)>0, s\neq 1,
 \end{array}
  \right.
\end{eqnarray*}
which can be continued
meromorphically to the whole complex $s$-plane, except for a simple pole at $s=1$ with its residue 1, see  Srivastava (2003) and  Choi et al. (2004) for details. One of the most celebrated formulas, discovered by Euler in 1734, is the following formula for positive even  integers
\begin{equation}
\zeta(2n)=(-1)^{n+1}\frac{2^{2n-1}}{(2n)!}\pi^{2n}B_{2n}, n\in\Bbb{N}_0,
\end{equation}
where $\Bbb{N}_0= \Bbb{N} \cup \{0\}$,   $B_n$ is  the   $n$th Bernoulli numbers. Here $\Bbb{N}$ is the set of natural
numbers. Since then,
many new proofs have  been obtained, see, for example, Titchmarsh and  Heath-Brown (1986),
 Amo et al. (2011),   Arakawa et al. (2014) and  Ribeiro (2018). On the contrary, however, no analogous closed forms
representation of $\zeta(s)$ at odd integers or   fractional points are known (cf. Srivastava and Choi (2012), P. 167). Even up to now, for positive odd integer arguments the Riemann zeta function can
only be expressed by series and integral. One possible
integral expression is established by
Cvijovi\'c and Klinowski (2002)  as follows
\begin{equation}
\zeta(2n+1)=(-1)^{n+1}\frac{(2\pi)^{2n+1}}{2(2n+1)!}\int_0^1 B_{2n+1}(u)\cot(\pi u)du, n\in\Bbb{N},
\end{equation}
 where $B_n(x)$   are  Bernoulli   polynomials  defined by the generating function (cf. Lu (2011))
$$\frac{te^{tx}}{e^t-1}=\sum_{n=0}^{\infty}B_n(x)\frac{t^n}{n!}, \; |t|<2\pi.$$
The Bernoulli numbers $B_n=B_n(0)$ are well-tabulated (see, for example, Srivastava(2003)):
$$B_0=1,B_1=-\frac12, B_2=\frac16, B_4=-\frac{1}{30}, B_6=\frac{1}{42},B_{2n+1}=0 \;(n=1,2,\cdots),\cdots,$$
from which one can finds  that
\begin{eqnarray*}
\zeta(2)=\frac{\pi^{2}}{6}, \zeta(4)=\frac{\pi^{4}}{90}, \zeta(6)=\frac{\pi^{6}}{945}, \zeta(8)=\frac{\pi^{8}}{9450},\cdots.
\end{eqnarray*}


The zeta  function $\zeta(s)$   has also the following integral
representation (cf.   Srivastava and  Choi (2012,  p.172))
\begin{eqnarray}
\zeta(s)=\frac{(1-2^{1-s})^{-1}}{\Gamma(s+1)}\int_0^{\infty}\frac{t^{s}e^{t}}{(e^t+1)^2}dt,\;\; {\cal{R}}(s)>0.
\end{eqnarray}
Note that  there is an extra 2 in (51) of  Srivastava and  Choi (2012, p.172).

 The aim of this note is to present a simple proof of   the  recurrence formula (1.1) for $\zeta(2n)$ and  the    integral representations for $\zeta(n)$ and $\zeta(n-\frac12)$ by making use of  probabilistic method.


\section{ The main results and their proofs}

In this section we present a new elementary proof to  the  following well known results.

\begin{proposition} For  Riemann's zeta function $\zeta$, we have
\begin{eqnarray}
\zeta(2n)=(-1)^{n-1}\frac{2^{2n-1}}{(2n)!}\pi^{2n}B_{2n},  n\in\Bbb{N}_0,
\end{eqnarray}
\begin{eqnarray}
\zeta\left(n-\frac{1}{2}\right)=\frac{2^{n}\int_{1}^{\infty}\frac{(\ln y)^{n-\frac{1}{2}}}{(1+y)^2}dy}{\sqrt{\pi}(2n-1)!!(1-2^{-\frac{2n-3}{2}})},
\;\; n\in \mathbb{N},
\end{eqnarray}
and
\begin{eqnarray}
\zeta\left(n\right)=\frac{(1-2^{1-n})^{-1}}{n!}\int_{0}^{\infty}\frac{x^n e^{-x}}{(1+e^{-x})^2}dx,
\;\; n\in \mathbb{N}, n>1,
\end{eqnarray}
where $B_n$ is  the   $n$th Bernoulli numbers.
\end{proposition}

To prove the  proposition, we need the following three lemmas.

\begin{lemma}
We assume that random variable $X$ has the standard logistic distribution with the  probability density function (pdf)
\begin{eqnarray}
f(x)=\frac{\exp(-x)}{(1+\exp(-x))^{2}}, -\infty<x<\infty.
\end{eqnarray}
 Then the  moment generating function (mgf) of $X$  is given by
\begin{eqnarray}
E[\exp(tX)]=\frac{\pi t}{\sin(\pi t)},\; |t|<1.
\end{eqnarray}
\end{lemma}

{\bf Proof}\;  By the definition of mgf, for any $|t|<1$,  we get
\begin{eqnarray*}
\begin{split}
E(e^{tX})&=\int_{-\infty}^{\infty} e^{tx}\frac{e^{-x}}{(1+e^{-x})^2}dx\\
&=\int_0^1 y^{-t}(1-y)^tdy=B(1-t,1+t)\\
&=\Gamma(1-t)\Gamma(1+t)\\
&=\frac{\pi t}{\sin(\pi t)},
\end{split}
\end{eqnarray*}
where $B(\cdot,\cdot)$ is the Beta function and $\Gamma(\cdot)$ is the $\Gamma$ function, we have used the fact that (see Gradshteyn and Ryzhik (1980) P.896))
$$\Gamma(1-t)\Gamma(t)=\frac{\pi}{\sin \pi t}$$ in the last  equality.

\begin{lemma}
We assume that random variable $X$ has the standard $1$-dimensional   elliptically symmetric logistic distribution with pdf
\begin{eqnarray}
f(x)=c\frac{\exp(-x^{2})}{(1+\exp(-x^{2}))^{2}}, -\infty<x<\infty,
\end{eqnarray}
where $$c=\left(\int_{0}^{\infty}t^{-\frac{1}{2}}\frac{e^{-t}}{(1+e^{-t})^{2}}dt\right)^{-1}.$$
 Then the characteristic function of $X$  is given by
\begin{eqnarray}
E[\exp(itX)]=1+\sum\limits_{n=1}^\infty(-1)^{n} \frac{c\sqrt{\pi}}{2^{2n+1}}\frac{t^{2n}}{n!}\left(1-2^{-\frac{2n-3}{2}}\right)\zeta\left(n-\frac{1}{2}\right),
\end{eqnarray}
where    $\zeta$ is the   Riemann zeta function.
\end{lemma}

{\bf Proof}\;   Let $h(x)=(1+\exp(-x^{2}))^{2}$, with the expansion for $x\neq 0$,
$$(1+\exp(-x^{2}))^{-2}=\sum\limits_{k=1}^\infty(-1)^{k-1}k\exp(-(k-1)x^{2}),$$
 we rewrite (2.6) as $$f(x)=c\sum\limits_{k=1}^\infty(-1)^{k-1}k\exp(-kx^{2}).$$
Noting that $f(-x)=f(x),-\infty<x<\infty$, so that all the odd-order moments of $f$ are zero.
Hence, we only need to determine the even-order moments. We get
\begin{eqnarray*}
\begin{split}
E(X^{2m}) & = 2\int_{0}^{\infty}x^{2m}f(x)dx \\
& = 2c\int_{0}^{\infty}x^{2m}\sum\limits_{k=1}^\infty(-1)^{k-1}k\exp(-kx^{2})dx \\
& = c\sum\limits_{k=1}^\infty(-1)^{k-1}\frac{\sqrt{\pi}}{2^{2m}}\frac{(2m)!}{m!}k^{-\frac{2m-1}{2}} \\
& = -c\sum\limits_{k=1}^\infty \frac{\sqrt{\pi}}{2^{2m}}\frac{(2m)!}{m!}(2k)^{-\frac{2m-1}{2}}+c\sum\limits_{k=1}^\infty \frac{\sqrt{\pi}}{2^{2m}}\frac{(2m)!}{m!}(2k-1)^{-\frac{2m-1}{2}} \\
& =c \sum\limits_{k=1}^\infty \frac{\sqrt{\pi}}{2^{2m}}\frac{(2m)!}{m!}k^{-\frac{2m-1}{2}}-c\sum\limits_{k=1}^\infty \frac{\sqrt{\pi}}{2^{2m}}2\frac{(2m)!}{m!}!(2k)^{-\frac{2m-1}{2}} \\
& = \frac{\sqrt{\pi}c}{2^{2m}}\frac{(2m)!}{m!}(1-2^{-\frac{2m-3}{2}})\zeta(m-\frac{1}{2}),
\end{split}
\end{eqnarray*}
where we have used the fact that $$\int_{0}^{\infty}\exp(-bx^{2})x^{2k}dx=\frac{\sqrt{\pi}}{2}\frac{1}{2}\frac{3}{2}\ldots\frac{2k-1}{2}b^{-\frac{2k+1}{2}}.$$
For any $t\in (-\infty,\infty)$,   we get the  characteristic function of $X$ by performing the following calculations
\begin{eqnarray*}
\begin{split}
E[\exp(itX)]
&=\int_{-\infty}^{\infty}\exp(itx)f(x)dx\\
&=E[\exp(itX)+\exp(-itX)]/2\\
&=E\left[1+\sum\limits_{n=1}^\infty(-1)^{n}\frac{t^{2n}X^{2n}}{(2n)!}\right]\\
&=1+\sum\limits_{n=1}^\infty(-1)^{n}\frac{t^{2n}E(X^{2n})}{(2n)!}\\
&=1+\sum\limits_{n=1}^\infty(-1)^{n} \frac{c\sqrt{\pi}}{2^{2n}}\frac{(2n)!}{n!}\frac{t^{2n}}{(2n)!}\left(1-2^{-\frac{2n-3}{2}}\right)\zeta\left(n-\frac{1}{2}\right)\\
&=1+\sum\limits_{n=1}^\infty(-1)^{n} \frac{c\sqrt{\pi}}{2^{2n}}\frac{t^{2n}}{n!}\left(1-2^{-\frac{2n-3}{2}}\right)\zeta\left(n-\frac{1}{2}\right).
\end{split}
\end{eqnarray*}
This ends the proof of Lemma 2.2.
\begin{lemma}
We assume that random variable $X$ has the standard half-logistic distribution with the  pdf
\begin{eqnarray}
f(x)=\frac{2\exp(-x)}{(1+\exp(-x))^{2}},\;\; x>0.
\end{eqnarray}
 Then the mgf of $X$  is given by
\begin{eqnarray}
E[\exp(tX)]=1+2\sum\limits_{n=1}^\infty (1-2^{1-n})\zeta(n)t^n,\;\; |t|<1,
\end{eqnarray}
where    $\zeta$ is the   Riemann zeta function.
\end{lemma}
{\bf Proof}\; The mean of $X$ is given by
$$E(X)=2\int_0^{\infty}\frac{xe^{-x}}{(1+\exp(-x))^{2}}dx=2\ln 2.$$
 Using the  expansion
\begin{eqnarray*}
f(x)&=&2\sum_{k=1}^{\infty}(-1)^{k-1}k e^{-kx}\\
&=&2\sum_{k=1}^{\infty}(2k-1) e^{-(2k-1)x}-2\sum_{k=1}^{\infty}2k e^{-2kx}, \; x>0,
\end{eqnarray*}
we get, for any positive integer $n>1$,
\begin{eqnarray*}
E(X^{n}) & = &2\int_{0}^{\infty}x^{n}f(x)dx \\
& =& 2  \sum_{k=1}^{\infty}(2k-1)\int_0^{\infty} x^n  e^{-(2k-1)x}dx- 2  \sum_{k=1}^{\infty}2k\int_0^{\infty} x^n  e^{-2kx}dx \\
&=&2n!  \sum_{k=1}^{\infty}\frac{1}{(2k-1)^n}-2n!  \sum_{k=1}^{\infty}\frac{1}{(2k)^n}\\
&=& 2n!(1-2^{1-n})\zeta(n).
\end{eqnarray*}
Then we have
\begin{eqnarray*}
\begin{split}
E(e^{tX})&=1+\sum\limits_{k=1}^\infty\frac{E(X^{k})}{k!}t^{k}\\
&=1+2t\ln 2+\sum\limits_{k=2}^\infty\frac{2k!\zeta(k)(1-2^{1-k})}{k!}t^k\\
&=1+2t\ln 2+2\sum\limits_{k=2}^{\infty} (1-2^{1-k})\zeta(k)t^k, \; |t|<1,
\end{split}
\end{eqnarray*}
where we have used the fact
$$\lim_{s\to 1}(s-1)\zeta(s)=1.$$
This completes the proof of Lemma 2.3.

{\bf Proof of  Proposition  2.1}

The  mgf of  the standard logistic distribution can  be written as
\begin{eqnarray}
E(e^{tX})=1+\sum\limits_{n=1}^\infty\frac{[2^{2n-1}-1]\zeta(2n)}{2^{2(n-1)}}t^{2n},
\end{eqnarray}
see, for example, Ghosh, Choi and  Li (2010).   Comparing   (2.4) and (2.7) yields  $g(t)=h(t), |t|<1$, where
$$g(t)=1+\sum\limits_{n=1}^\infty\frac{[2^{2n-1}-1]\zeta(2n)}{2^{2(n-1)}}t^{2n},$$
and $$h(t)=\frac{\pi t}{\sin(\pi t)}.$$
 Using the series expansion
$$\frac{\pi t}{\sin \pi t}=\sum\limits_{k=0}^\infty(-1)^{k-1}\frac{2^{2k}-2}{(2k)!}B_{2k}(\pi t)^{2k},$$
where   $B_{2k}$ is  the   $2k$th Bernoulli numbers,
we have
$$\sum\limits_{n=1}^\infty\frac{[2^{2n-1}-1]\zeta(2n)}{2^{2(n-1)}}t^{2n}= \sum\limits_{k=1}^\infty(-1)^{k-1}\frac{2^{2k}-2}{(2k)!}B_{2k}(\pi t)^{2k},\;\; |t|<1,$$
from which we deduce that
 \begin{eqnarray*}
\zeta(2n)=(-1)^{n-1}\frac{2^{2n-1}}{(2n)!}\pi^{2n}B_{2n}.
\end{eqnarray*}
This completes the proof of (2.1).

 Now we prove (2.2). Denoting by $$H(t)=\int_{-\infty}^{\infty}\exp(itx)f(x)dx,$$
 and  $$G(t)=1+\sum\limits_{k=1}^\infty(-1)^{k}\frac{c\sqrt{\pi}}{2^{2k+1}}\frac{t^{2k}}{k!}\left(1-2^{-\frac{2k-3}{2}}\right)\zeta\left(k-\frac{1}{2}\right),$$
 where $f$ is defined by (2.4).
 Taking $2n$th and $(2n+1)$th derivatives of the two functions with respect to $t$, we get
$$H^{(2n)}(t)=(-1)^{n}2c\int_{0}^{\infty}x^{2n}\cos tx\frac{\exp(-x^{2})}{(1+\exp(-x^{2}))^{2}}dx,$$
$$H^{(2n+1)}(t)=(-1)^{n+1}2c\int_{0}^{\infty}x^{2n-1}\sin tx\frac{\exp(-x^{2})}{(1+\exp(-x^{2}))^{2}}dx,$$
and
$$G^{(2n-1)}(t)=\sum\limits_{k=n}^\infty(-1)^{k}
\frac{c\sqrt{\pi}\prod_{l=0}^{2n-2}(2k-l)}{2^{2k}k!}\left(1-2^{-\frac{2k-3}{2}}\right)\zeta\left(k-\frac{1}{2}\right)t^{2k-2n+1},$$
$$G^{(2n)}(t)=\sum\limits_{k=n}^\infty(-1)^{k}
\frac{c\sqrt{\pi}\prod_{l=0}^{2n-1}(2k-l)}{2^{2k}k!}\left(1-2^{-\frac{2k-3}{2}}\right)\zeta\left(k-\frac{1}{2}\right)t^{2k-2n}.$$
Note that  $H(t)=G(t)$ for any real $t$, and thus  $H^{(n)}(t)=G^{(n)}(t)$ for any real $t$ and any positive integers $n$. In particular, $H^{(n)}(0)=G^{(n)}(0)$. However, $H^{(2n+1)}(0)=G^{(2n+1)}(0)=0$,   and  from  $H^{(2n)}(0)=G^{(2n)}(0)$ we have
\begin{eqnarray*}
\begin{aligned}
 \zeta\left(\frac{2n-1}{2}\right)&=\frac{2^{2n+1}n!\int_{0}^{\infty}x^{2n}\frac{\exp(-x^2)}{(1+\exp(-x^2))^{2}}dx}{\sqrt{\pi}(2n)!(1-2^{-\frac{2n-3}{2}})}\\
 &=\frac{2^{n}\int_{0}^{\infty}x^{\frac{2n-1}{2}}\frac{\exp(-x)}{(1+\exp(-x))^{2}}dx}{\sqrt{\pi}(2n-1)!!(1-2^{-\frac{2n-3}{2}})}\\
 &= \frac{2^{n}\int_{1}^{\infty}\frac{(\ln y)^{n-\frac{1}{2}}}{(1+y)^2}dy}{\sqrt{\pi}(2n-1)!!(1-2^{-\frac{2n-3}{2}})}, \;\; n\in \mathbb{N},
\end{aligned}
\end{eqnarray*}
 which concludes the proof of (2.2).

 Finally we prove (2.3). Using (2.9) one has
 \begin{equation}
 2\int_0^{\infty}e^{tx}\frac{e^{-x}}{(1+e^{-x})^2}dx=1+2t\ln 2+2\sum\limits_{k=2}^{\infty} (1-2^{1-k})\zeta(k)t^k, \; |t|<1.
 \end{equation}
Taking $n$th  derivative of both sides of (2.11)  with respect to $t$ and then setting $t=0$ yields the desired result.

\noindent{\bf Acknowledgements.}
  The research was supported by the National Natural Science Foundation of China (No. 11571198).

\end{document}